\documentstyle[twoside]{article}
\begin{document}
\def\f#1#2{\mbox{${\textstyle \frac{#1}{#2}}$}}
\newcommand{\fs}{\f{1}{2}}
\begin{center}

{\bf \Large On an Extension of Extended Beta and Hypergeometric Functions}\\[9mm]

{\bf R. K. Parmar$^{1}$, 
P. Chopra$^{2}$ and R. B. Paris$^{3}$}\\[3mm]

$^{1}$Department of Mathematics, Government
College of Engineering and Technology, \\
Bikaner 334004, Rajasthan State, India\\

E-Mail: rakeshparmar27@gmail.com\\[2mm]

$^{2}$Department of Mathematics, Marudhar Engineering College,\\
Bikaner-334001, Rajasthan State, India\\

E-Mail: purnimachopra@rediffmail.com\\[2mm]

$^{3}$School of Engineering, Computing and Applied Mathematics,\\
 University of Abertay Dundee, Dundee DD1 1HG, UK\\

E-Mail: r.paris@abertay.ac.uk\\

\bigskip

\textbf{Abstract}
\end{center}
Motivated mainly by certain interesting recent extensions of the Gamma, Beta and hypergeometric functions, we introduce here new extensions of the Beta function, hypergeometric and confluent hypergeometric functions. We systematically investigate several properties of each of these extended functions, namely their various integral representations, Mellin transforms, derivatives, transformations, summation formulas and asymptotics. Relevant connections of certain special cases of the main results presented here are also pointed out.

\bigskip

\noindent
{\bf 2010 Mathematics Subject Classification.}  Primary 33B20, 33C20; Secondary 33B15, 33C05.\\

\noindent
{\bf Key Words and Phrases.} Generalized Gamma function; Extended Gamma function; Extended Beta function; Extended hypergeometric function; Extended confluent hypergeometric function; Mellin transforms; Transformation formulas; Summation formula.
\vspace{0.6cm}

\begin{center}
{\bf 1. \ Introduction, Definitions and Preliminaries}
\end{center}
\setcounter{section}{1}
\setcounter{equation}{0}
\renewcommand{\theequation}{\arabic{section}.\arabic{equation}}

\noindent 
Extensions of a number of well-known special functions have been investigated recently by several authors (see \cite{Ch-Qa-Ra-Zu, Ch-Qa-Sr-Pa, Ch-Te-Ve, Ch-Zu-94, Ch-Zu-95,
 Ch-Zu-02, Ch-Zu-02-Book}). In 1994, Chaudhry and Zubair \cite{Ch-Zu-94} introduced the following extensions of the incomplete Gamma functions in the form

\begin{equation}\label{Incomplete-gamma-p}
\gamma(\alpha,x;p):=\int_0^x\, t^{\alpha-1}\,\exp \left(-t- \frac{p}{t}\right)dt\qquad(\Re(p)>0; p=0, \ \Re(\alpha)>0)
\end{equation}
and
\begin{equation}\label{Incomplete-Gamma-p}
\Gamma(\alpha,x;p):=\int_x^\infty\, t^{\alpha-1}\,\exp \left(-t- \frac{p}{t}\right)dt\qquad(\Re(p)\geq 0)
\end{equation}
with $|\arg\,x|<\pi$.
These functions satisfy the following decomposition formula:
\begin{eqnarray}
\gamma (\alpha,x;p) + \Gamma (\alpha,x;p)&\equiv&\Gamma_p (x)=\int_0^\infty\, t^{\alpha-1}\, \exp \left(-t- \frac{p}{t}\right)dt\nonumber\\
&=&2p^{\alpha/2}K_{\alpha}\left(2\sqrt{p}\right) \qquad (\Re (p)>0),\label{gamma+Gamma-function-p}
\end{eqnarray}
where, here and throughout, $K_\nu(x)$ denotes the modified Bessel function.
In 1997, Chaudhry {\it et al.} \cite[Eq.~(1.7)]{Ch-Qa-Ra-Zu} presented the following extension of the Beta function given by
\begin{equation}\label{Beta-function-p}
 B(x,\,y\,;\,p):=\int_0^1\, t^{x-1}\,(1-t)^{y-1}\, \exp \left(- \frac{p}{t(1-t)}\right)\,dt
\end{equation}
$$(\Re (p)>0;\, p=0,\,\Re (x)>0,\,\Re (y)>0)$$
and showed that this extension has certain connections with the Macdonald, error  and Whittaker functions.

More recently, Chaudhry {\it et al.} \cite{Ch-Qa-Sr-Pa}
used $B(x,\,y\,;\,p)$ to extend the hypergeometric and the confluent hypergeometric functions as follows:
\begin{equation}\label{GHF-p}
 F_p\,(a,\, b;\,c;\, z):=\sum_{n=0}^\infty\,(a)_n\, \frac{B(b+n,\,c-b\,;\,p)}{B(b,\,c-b)}\,\frac{z^n}{n!}
\end{equation}
  $$ \left(p \geq 0, \,\, |z|<1\,;\, \Re(c)>\Re(b)>0  \right)$$
and
\begin{equation}\label{CHF-p}
 \Phi_p\,(b;\,c;\, z):=\sum_{n=0}^\infty\, \frac{B(b+n,\,c-b\,;\,p)}{B(b,\,c-b)}\,\frac{z^n}{n!}
\end{equation}
  $$ \left(p \geq 0 \,;\, \Re(c)>\Re(b)>0  \right).$$
Among several interesting and potentially useful properties of the extended hypergeometric functions defined by (\ref{GHF-p}) and (\ref{CHF-p}), the following integral representations were also given by Chaudhry {\it et al.} \cite[Eq.~(3.2)]{Ch-Qa-Sr-Pa} and \cite[Eq.~(3.6)]{Ch-Qa-Sr-Pa}:

\begin{equation}\label{GHF-Integral-p}
 F_{p}\,(a,\, b;\,c;\, z):=\frac{1}{B(b,\,c-b)}\,\int_0^1\,t^{b-1}\,(1-t)^{c-b-1}\, (1-zt)^{-a}\, \exp \left(- \frac{p}{t(1-t)}\right)\,dt
\end{equation}
 $$ \left(|\arg (1-z)|<\pi;\,\,p > 0;\,p=0\,\,\mbox{and}\,\, \, \Re(c)>\Re(b)>0  \right)$$
and
\begin{equation}\label{CHF-Integral-p}
 \Phi_{p}\,(b;\,c;\, z):=\frac{1}{B(b,\,c-b)}\,\int_0^1\,t^{b-1}\,(1-t)^{c-b-1}\, \exp \left(zt- \frac{p}{t(1-t)}\right)\,dt
\end{equation}
 $$ \left(p > 0;\,p=0\,\,\mbox{and}\,\,  \Re(c)>\Re(b)>0  \right),$$
respectively.
They investigated these functions and gave their various
integral representations, beta distribution, certain properties including differentiation formulas, Mellin transform, transformation formulas,
recurrence relations, summation formula, asymptotic formulas and certain interesting connections with some 
well-known special functions.
It is clearly seen, that in the case $p=0$, (\ref{gamma+Gamma-function-p}), (\ref{Beta-function-p}), (\ref{GHF-p}) and (\ref{CHF-p}) reduce to the usual Gamma function, Beta function, Gauss hypergeometric and confluent
hypergeometric functions, respectively.

In 1997, Chaudhry and Zubair \cite{Ch-Zu-97} considered the following extensions of the incomplete gamma functions (\ref{Incomplete-gamma-p}) and (\ref{Incomplete-Gamma-p}) in the form
\begin{equation}\label{Incomplete-gamma-K-p}
\gamma_{\nu}(\alpha,x;p):=\sqrt{\frac{2p}{\pi}}\,\int_0^x\, t^{\alpha-\frac{3}{2}}\, \exp \left(-t\right)K_{\nu+\frac{1}{2}}\left(\frac{p}{t}\right) \,dt
\end{equation}
and
\begin{equation}\label{Incomplete-Gamma-K-p}
\Gamma_{\nu}(\alpha,x;p):=\sqrt{\frac{2p}{\pi}}\,\int_x^\infty\, t^{\alpha-\frac{3}{2}}\, \exp \left(-t\right)K_{\nu+\frac{1}{2}}\left(\frac{p}{t}\right) \,dt
\end{equation}
$$(|\arg\,x|<\pi,\,\Re (p)>0,-\infty<\alpha<\infty).$$
respectively. These latter functions satisfy the following decomposition formula:
\[\gamma_{\nu}(\alpha,x;p)+\Gamma_{\nu}(\alpha,x;p)=2^{\alpha-2}\pi^{-1}\sqrt{p}\hspace{6cm}\]
\begin{equation}\label{gamma+Gamma-function-K-p}
\hspace{15mm} \times G_{4,0}^{0,4}\left(\frac{p^{2}}{16}\mid \frac{1}{2}\left(\nu+\frac{1}{2}\right),-\frac{1}{2}\left(\nu+\frac{1}{2}\right),\frac{1}{2}\left(\alpha+\frac{1}{2}\right),\frac{1}{2}\left(\alpha-\frac{1}{2}\right)\right)
\end{equation}
$$(\Re (p)>0,-\infty<\alpha<\infty),$$
where $G$ is the Meijer $G$-function.
From the fact that \cite[Eq.~(10.39.2)]{DLMF}
\begin{equation}\label{K-case}
K_{\frac{1}{2}}\left(z\right)=\sqrt{\frac{\pi}{2z}}\,e^{-z},
\end{equation}
it is easily seen that (\ref{Incomplete-gamma-K-p}) and (\ref{Incomplete-Gamma-K-p}) reduce to the generalized incomplete gamma functions (\ref{Incomplete-gamma-p}) and (\ref{Incomplete-Gamma-p}) when $\nu=0$.

Motivated essentially by the demonstrated potential for applications of these generalized incomplete gamma functions in many diverse areas of mathematical, physical, engineering and statistical sciences (see, for details, \cite{Ch-Zu-97,Ch-Zu-98} and the references cited therein),
we introduce here another interesting extension of the extended Beta function $B(x,\,y\,;\,p)$ given by
\begin{equation}\label{Beta-function-K}
B_{\nu}(x,\,y\,;\,p):= \sqrt{\frac{2p}{\pi}}\,\int_{0}^{1}\, t^{x-\frac{3}{2}}\,(1-t)^{y-\frac{3}{2}}\, K_{\nu+\frac{1}{2}}\left(\frac{p}{t(1-t)}\right) \,dt \qquad (\Re (p)>0)
\end{equation}
We note that this extension preserves the symmetry in the parameters $x$ and $y$, and we have
$$B_{\nu}(x,\,y\,;\,p)=B_{\nu}(y,\,x\,;\,p).$$

\vspace{2mm}

\noindent
\textbf{Remark 1.} The special case of (\ref{Beta-function-K}) corresponding to $\nu=0$  
is easily seen to reduce to the extended beta function (\ref{Beta-function-p}) upon making use of (\ref{K-case}).
\vspace{0.3cm}

Further, making use of the extended Beta function $B_{\nu}(x,\,y\,;\,p)$ in (\ref{Beta-function-K}), we consider other extensions of the extended Gauss hypergeometric and the confluent hypergeometric functions. For each of these new extensions we obtain various integral representations, properties and Mellin transforms, together with differentiation, transformation, summation and asymptotic formulas. Relevant connections of certain special cases of the main results presented here are also pointed out.

\vspace{0.6cm}

\begin{center}
{\bf 2. \ Integral Representations of $B_{\nu}(x,\,y\,;\,p)$} 
\end{center}
\setcounter{section}{2}
\setcounter{equation}{0}
\renewcommand{\theequation}{\arabic{section}.\arabic{equation}}

\noindent
In this section, we obtain integral representations, a representation of $B_{\nu}(x,\,y\,;\,p)$ for integer values of $\nu$ and certain special cases.

\vspace{3mm}

\noindent
\textbf{Theorem 1.} \emph{The following various integral representations for} $B_{\nu}(x,\,y\,;\,p)$ in (\ref{Beta-function-K}) \emph{hold true}:
\begin{equation}\label{Beta-Integral-1}
B_{\nu}(x,\,y\,;\,p)= 2\sqrt{\frac{2p}{\pi}}\,\int_0^{\frac{\pi}{2}} \, \cos^{2(x-1)}\theta\, \sin^{2(y-1)}\theta K_{\nu+\frac{1}{2}}\left(p\,\sec^2\theta\,\csc^2 \theta\right)d\theta,
\end{equation}

\begin{equation}\label{Beta-Integral-2}
B_{\nu}(x,\,y\,;\,p)=\sqrt{\frac{2p}{\pi}}\,\int_0^\infty\, \frac{u^{x-\frac{3}{2}}}{(1+u)^{x+y-1}}\, K_{\nu+\frac{1}{2}}\left(p(2+u+\frac{1}{u})\right)du
\end{equation}
{\it and}
\begin{equation}\label{Beta-Integral-3}
B_{\nu}(x,\,y\,;\,p)=  2^{2-x-y}\,\sqrt{\frac{2p}{\pi}}\, \int_{-1}^1\, (1+u)^{x-\frac{3}{2}}\,(1-u)^{y-\frac{3}{2}} K_{\nu+\frac{1}{2}}\left(\frac{4p}{1-u^2}\right)du
\end{equation}
$$ \left(\Re(p)>0\right).$$

\noindent{\it Proof.}\ \ 
Equations (\ref{Beta-Integral-1}), (\ref{Beta-Integral-2}) and (\ref{Beta-Integral-3}) can be obtained by employing the transformations
$t= \cos^2 \theta$, $t= u/(1+u)$ and $t=(1+u)/2$ in (\ref{Beta-function-K}), respectively.
\hfill$\Box$

\vskip 3mm

\noindent
\textbf{Remark 2.} Clearly, when $\nu=0$, Theorem 1 reduces to the corresponding integrals in \cite{Ch-Qa-Ra-Zu}.

\vskip 3mm

\noindent
\textbf{Theorem 2.} \emph{The following representation for integer values of} $\nu=n$ in (\ref{Beta-function-K})
\emph{holds true}:
\begin{equation}\label{Beta-Integral-values}
B_{n}(x,\,y\,;\,p)= \sum_{m=0}^{n}\frac{(2p)^{-m}}{m!}\;\frac{\Gamma (n+m+1)}{\Gamma(n-m+1)}\,B(x+m,\,y+m\,;\,p)
\end{equation}
$$(n \in N\,\,\mbox{and}\,\,\Re(p)>0 ).$$

\noindent{\it Proof.}\ \  Using the fact (see, for example, \cite[Eq.~(10.49.12)]{DLMF})
\begin{equation}\label{K-formula}
 K_{n+\frac{1}{2}}\left(z\right)=\sqrt{\frac{\pi}{2z}}\,e^{-z}\sum_{m=0}^{n}\frac{(2z)^{-m}}{m!} \;\frac{\Gamma (n+m+1)}{\Gamma (n-m+1)}
\end{equation}
in (\ref{Beta-function-K}) and applying the definition (\ref{Beta-function-p}), we obtain the desired representation (\ref{Beta-Integral-values}). \hfill $\Box$

\vspace{3mm}

\noindent
\textbf{Remark 3.}  The Meijer $G$-\emph{function} \cite[p.~232, Eq. (5.124)]{Ch-Zu-02-Book}, the Whittaker function $W_{\alpha,\beta}(z)$ \cite[p.~334]{DLMF},  and the confluent function $U(a,b,z)$ \cite[p.~325]{DLMF} are expressible in terms of the extended beta function $B(x,\,y\,;\,p)$ for $\Re(p)>0$ as follows:

\begin{equation}\label{G-function}
B(x,\,y\,;\,p)=\sqrt{\pi}\,2^{1-x-y}\,G_{2,3}^{3,0} \left(4p \mid
\begin{array}{rrr}
\frac{x+y}{2}\,,\frac{x+y+1}{2} \\
 \,\, 0,\,x,\,y
\end{array}
\right)
\end{equation}
\begin{equation}\label{Whittaker-function}
B(x,\,x\,;\,p) = \sqrt{\pi}\,2^{-x}p^{\frac{x-1}{2}}\,e^{-2p }\,W_{-\frac{x}{2},\frac{x}{2}}(4p)
\end{equation}
and
\begin{equation}\label{Tricomi-function}
B(x,\,x\,;\,p)= \sqrt{\pi}\,2^{1-2x}\,e^{-4p}\, U({\fs,1-x},4p).
\end{equation}

Now, applying the relationships (\ref{G-function}) -- (\ref{Tricomi-function}) to (\ref{Beta-Integral-values}),
we can deduce certain interesting representations of the extended beta function in
(\ref{Beta-function-K}) for integer values of $\nu$. These are given in Corollary 1 below, where we state the resulting representations without proof.

\vskip 3mm

\noindent
\textbf{Corollary 1.} \emph{Each of the following representations for integer values of} $\nu=n$ in (\ref{Beta-Integral-values}) \emph{holds true}:
\[B_{n}(x,\,y\,;\,p)= \sqrt{\pi}\,2^{1-x-y}\,\sum_{m=0}^{n} \frac{(8p)^{-m}}{m!} \;\frac{\Gamma (n+m+1) }{\Gamma (n-m+1)}\hspace{4cm}\]
\begin{equation}\label{Beta-Integral-values-0}
\hspace{3cm}  \times G_{2,3}^{3,0} \left(4p \mid
\begin{array}{rrr}
\frac{x+y+2m}{2}\,,\frac{x+y+2m+1}{2} \\
 \,\, 0,\,x+m\,,y+m
\end{array}
\right),
\end{equation}
\begin{equation}\label{Beta-Integral-values-1}
B_{n}(x,\,x\,;\,p)= \sqrt{\pi}\,2^{-x}p^{\frac{x-1}{2}}\,e^{-2p }\,\sum_{m=0}^{n} \frac{(4p\sqrt{p})^{-m}}{m!} \;\frac{\Gamma (n+m+1) }{\Gamma (n-m+1)}\, W_{-\frac{x+m}{2},\frac{x+m}{2}}(4p)
\end{equation}
{\it and}
\begin{equation}\label{Beta-Integral-values-2}
B_{n}(x,\,x\,;\,p)= \sqrt{\pi}\,2^{1-2x}\,e^{-4p}\,\sum_{m=0}^{n} \frac{(8p)^{-m}}{m!} \;\frac{\Gamma (n+m+1) }{\Gamma (n-m+1)}\, U({\fs,1-x-m},4p).
\end{equation}
\vspace{0.6cm}

\begin{center}
{\bf 3. \ Some Properties of $B_{\nu}(x,\,y\,;\,p)$} 
\end{center}
\setcounter{section}{3}
\setcounter{equation}{0}
\renewcommand{\theequation}{\arabic{section}.\arabic{equation}}

\noindent
In this section, we obtain certain properties including, for example, functional relations and summation formulas for $B_{\nu}(x,\,y\,;\,p)$ as follows:

\vspace{3mm}

\noindent
\textbf{Theorem 3.} \emph{The following functional relation for} $B_{\nu}(x,\,y\,;\,p)$ in (\ref{Beta-function-K})
\emph{holds true}:
 \begin{equation}\label{functional-relation}
B_{\nu}(x+1,\,y\,;\,p) + B_{\nu}(x,\,y+1\,;\,p)=B_{\nu}(x,\,y\,;\,p).
\end{equation}

\noindent{\it Proof.}\ \ The left-hand side of (\ref{functional-relation}) becomes
$$B_{\nu}(x+1,\,y\,;\,p) + B_{\nu}(x,\,y+1\,;\,p)=\sqrt{\frac{2p}{\pi}}\,\int_0^1\, t^{x-\frac{3}{2}}\,(1-t)^{y-\frac{3}{2}}\left[t+(1-t)\right]K_{\nu+\frac{1}{2}}\left(\frac{p}{t(1-t)}\right) \,dt$$
which, after a simple algebraic manipulation and application of (\ref{Beta-function-K}), establishes the desired result. \hfill $\Box$.

\vspace{3mm}

\noindent
\textbf{Theorem 4.} \emph{The following summation formula for} $B_{\nu}(x,\,y\,;\,p)$ in (\ref{Beta-function-K})
\emph{holds true}:
 \begin{equation}\label{summation-formula-1}
B_{\nu}(x,\,1-y\,;\,p)= \sum_{n=0}^\infty\, \frac{(y)_n}{n!}\, B_{\nu}(x+n,\,1\,;\,p) \qquad (\Re(p)>0).
 \end{equation}
 
\noindent{\it Proof.}\ \  Application of the binomial theorem
$$ (1-t)^{-y} = \sum_{n=0}^\infty\, (y)_n\, \frac{t^n}{n!} \qquad (|t|<1)$$
to the definition (\ref{Beta-function-K}) of $B_{\nu}(x,\,y\,;\,p)$ yields
\[ B_{\nu}(x,\,1-y\,;\,p)= \sqrt{\frac{2p}{\pi}}\,\int_0^1\,\sum_{n=0}^\infty\,\frac{(y)_n}{n!}\,t^{x+n-\frac{3}{2}}(1-t)^{-\frac{1}{2}}\,
   K_{\nu+\frac{1}{2}}\left(\frac{p}{t(1-t)}\right)\,dt.\]
Interchange of the order of integration and summation in the last expression and use of (\ref{Beta-function-K}) then proves the desired identity. \hfill $\Box$

\vskip 3mm

\noindent
\textbf{Theorem 5.} \emph{The following infinite summation formula for} $B_{\nu}(x,\,y\,;\,p)$ in (\ref{Beta-function-K}) \emph{holds true}:
\begin{equation}\label{summation-formula-2}
B_{\nu}(x,\,y\,;\,p)= \sum_{n=0}^\infty\,  B_{\nu}(x+n,\,y+1\,;\,p) \qquad (\Re(p)>0).
\end{equation}

\noindent{\it Proof.}\ \  Replacing $(1-t)^{y-1}$ in (\ref{Beta-function-K}) by its series representation
$$ (1-t)^{y-1} = (1-t)^y\, \sum_{n=0}^\infty\,t^n\qquad (|t|<1),$$
we obtain
$$B_{\nu}(x,\,y\,;\,p)=\sqrt{\frac{2p}{\pi}} \int_0^1\, (1-t)^{y-\frac{1}{2}}\, \sum_{n=0}^\infty\,t^{x+n-\frac{3}{2}}\, K_{\nu+\frac{1}{2}}\left(\frac{p}{t(1-t)}\right)\,dt.$$
Interchange of the order of integration and summation in the last expression and use of (\ref{Beta-function-K})
then proves the desired identity. \hfill $\Box$
\vskip 3mm

\noindent
\textbf{Remark 4.}  In the special case $\nu=0$, (\ref{functional-relation}) -- (\ref{summation-formula-2})
reduce to the corresponding results given in \cite{Ch-Qa-Ra-Zu}.
\vspace{0.6cm}

\begin{center}
{\bf 4. \ Mellin Transforms of $B_{\nu}(x,\,y\,;\,p)$} 
\end{center}
\setcounter{section}{4}
\setcounter{equation}{0}
\renewcommand{\theequation}{\arabic{section}.\arabic{equation}}

\noindent
The Mellin transform  of a suitably integrable
function $f(t)$ with index $s$ is defined, as usual, by \cite{Er-Ma-Ob-Tr-I}
\begin{equation}\label{Mellin-transform}
{\mathcal M}\{f(\tau):\tau\to s\}:= \int_0^{\infty}\tau^{s-1} f(\tau)\, d\tau
\end{equation}
whenever the improper integral in (\ref{Mellin-transform}) exists.

\vspace{3mm}

\noindent
\textbf{Theorem 6.} \emph{The following Mellin transformation formula for} $B_{\nu}(x,\,y\,;\,p)$ in (\ref{Beta-function-K}) \emph{holds true}:
\begin{equation}\label{Mellin-transform-Beta}
{\mathcal M}\left\{B_{\nu}(x,\,y\,;\,p):p \to s\right\}=\frac{2^{s-1}}{\surd\pi}\,\Gamma(\frac{s-\nu}{2})\Gamma(\frac{s+\nu+1}{2})\;B(x+s,y+s)
\end{equation}
$$\big(\Re(x+s)>0,\; \Re(y+s)>0,\; \Re(s-\nu)>0,\; \Re(s+\nu)>-1\big),$$
where $B(x,y)$ denotes the classical Beta function,
provided that each member of (\ref{Mellin-transform-Beta}) exists.

\noindent{\it Proof.}\ \  Using the Definition (\ref{Mellin-transform}) of the
Mellin transform, we find from (\ref{Beta-function-K}) that
\[{\mathcal M}\left\{B_{\nu}(x,\,y\,;\,p):p \to s\right\}\hspace{8cm}\]
\begin{eqnarray}
&=&\int_0^{\infty}p^{s-1}
\left\{\sqrt{\frac{2p}{\pi}}\,\int_{0}^{1}\, t^{x-\frac{3}{2}}\,(1-t)^{y-\frac{3}{2}}\, K_{\nu+\frac{1}{2}}\left(\frac{p}{t(1-t)}\right) dt\right\}dp\nonumber\\
&=&\sqrt{\frac{2}{\pi}}\,\int_0^{\infty}p^{s-\frac{1}{2}}\left\{
\int_0^1 t^{x-\frac{3}{2}}\,(1-t)^{y-\frac{3}{2}}\, K_{\nu+\frac{1}{2}}\left(\frac{p}{t(1-t)}\right) dt\right\}dp\,.\label{Mellin-transform-Beta-1}
\end{eqnarray}
Upon interchanging the order of integration in (\ref{Mellin-transform-Beta-1}),
which can easily be justified by absolute convergence of the integrals
involved under the constraints stated in (\ref{Mellin-transform-Beta}), we find
\[{\mathcal M}\left\{B_{\nu}(x,\,y\,;\,p):p \to s\right\}\hspace{8cm}\]
\begin{eqnarray}
&=&\sqrt{\frac{2}{\pi}}\,\int_0^1t^{x-\frac{3}{2}}\,(1-t)^{y-\frac{3}{2}}
\left\{\int_0^{\infty}p^{s-\frac{1}{2}}K_{\nu+\frac{1}{2}}\left(\frac{p}{t(1-t)}\right)dp\right\} dt\nonumber \\
&=&\int_0^1t^{x+s-1}\,(1-t)^{y+s-1}
\left\{\sqrt{\frac{2}{\pi}}\,\int_0^{\infty}w^{s-\frac{1}{2}}K_{\nu+\frac{1}{2}}\left(w\right) dw\right\} dt,\label{Mellin-transform-Beta-2}
\end{eqnarray}
where we have set $w=p t/(1-t)$ in the inner $p$-integral.

The inner integral may be evaluated with the help of \cite[Eq.~(10.43.19)]{DLMF}
\begin{equation}\label{K-integral}
\int_0^{\infty}t^{s-\frac{1}{2}}K_{\alpha+\frac{1}{2}}\left(t\right) dt=2^{s-\frac{3}{2}}\,\Gamma(\frac{s-\alpha}{2})\Gamma(\frac{s+\alpha+1}{2})
\end{equation}
for $\Re(s-\alpha)>0$, $\Re(s+\alpha)>-1$,
and the integral with respect to $t$ can then be evaluated in terms of the classical Beta function (given by (\ref{Beta-function-p}) with $p=0$).
This evidently completes
our derivation of the Mellin transform formula (\ref{Mellin-transform-Beta}). \hfill $\Box$

\vspace{0.6cm}

\begin{center}
{\bf 5. \ The Beta Distribution of $B_{\nu}(x,\,y\,;\,p)$} 
\end{center}
\setcounter{section}{5}
\setcounter{equation}{0}
\renewcommand{\theequation}{\arabic{section}.\arabic{equation}}
\noindent
As a statistical application of $B_{\nu}(x,\,y\,;\,p)$, we define the extended beta distribution by employing $B_{\nu}(x,\,y\,;\,p)$, where the parameters $a$ and $b$ satisfy $\infty<a<\infty$,\,$\infty<b<\infty$ and $p>0$ as:
 \begin{equation}\label{Beta-Distribution}
  f(t)=\left\{ \begin{array}{ll} \displaystyle{ \frac{1}{B_{\nu}(a,\,b\,;\,p)} \sqrt{\frac{2p}{\pi}}\,t^{a-\frac{3}{2}}\,(1-t)^{b-\frac{3}{2}}\,
K_{\nu+\frac{1}{2}}\left(\frac{p}{t(1-t)}\right)} & (0<t<1),\\
\\
    0, & \mbox{otherwise.}
\end{array} \right.
\end{equation}
If $r$ is any real number, then we have the $r$th moment of $X$
   \begin{equation}\label{Expectation}
   E\left(X^r\right) = \frac{B_{\nu}(a+r,\,b\,;\,p)}{B_{\nu}(a,\,b\,;\,p)}
       \end{equation}
    $$ \left(p>0,\, -\infty <a<\infty,\, -\infty <b<\infty\right). $$
The particular case of (\ref{Expectation}) when $r=1$,
 \begin{equation}\label{Mean}
\mu =E\left(X\right)=\frac{B_{\nu}(a+1,\,b\,;\,p)}{B_{\nu}(a,\,b\,;\,p)},
    \end{equation}
represents the mean of the distribution, and
 \begin{equation}\label{Variance}
\sigma^2 =E\left(X^2\right)- \left\{E\left(X\right)\right\}^2
=\frac{B_{\nu}(a,\,b\,;\,p)B_{\nu}(a+2,\,b\,;\,p)- B_{\nu}^2(a+1,\,b\,;\,p)}{B_{\nu}^2(a,\,b\,;\,p)}
    \end{equation}
is the variance of the distribution.

The moment generating function of the distribution is
 \begin{equation}\label{MGF}
M(t)= \sum_{n=0}^\infty\, \frac{t^n}{n!}\, E\left(X^n\right)= \frac{1}{B_{\nu}(a,\,b\,;\,p)}\,
  \sum_{n=0}^\infty\,B_{\nu}(a+n,\,b\,;\,p)\,\frac{t^n}{n!}.
   \end{equation}
The cumulative distribution of (\ref{Beta-Distribution}) can be expressed as
 \begin{equation}\label{Cumulative-Distribution}
F(x)= \frac{B_{\nu,x}(a,\,b\,;\,p)}{B_{\nu}(a,\,b\,;\,p)},
   \end{equation}
where
 \begin{equation}\label{Incomplete-Beta-function}
B_{x,\nu}(a,\,b\,;\,p)=\sqrt{\frac{2p}{\pi}}\,\int_0^x\, t^{a-\frac{3}{2}}\,(1-t)^{b-\frac{3}{2}}\, K_{\nu+\frac{1}{2}}\left(\frac{p}{t(1-t)}\right)dt
   \end{equation}
$$ \left(p>0,\, -\infty <a<\infty,\, -\infty <b<\infty\right) $$
is the extended incomplete beta function.

\vspace{0.6cm}

\begin{center}
{\bf 6. \ Extension of Extended Hypergeometric Functions} 
\end{center}
\setcounter{section}{6}
\setcounter{equation}{0}
\renewcommand{\theequation}{\arabic{section}.\arabic{equation}}
\noindent

In this section, we extend the extended
Gauss hypergeometric and confluent hypergeometric functions in (\ref{GHF-p}) and (\ref{CHF-p})    by making use of $B_{\nu}(x,\,y\,;\,p)$ in (\ref{Beta-function-K}) as follows:
\begin{equation}\label{GHF-K}
 F_{p,\nu}\,(a,\, b;\,c;\, z):=\sum_{n=0}^\infty\,(a)_n\,\frac{B_{\nu}(b+n,\,c-b\,;\,p)}{B(b,\,c-b)}\,\frac{z^n}{n!}
\end{equation}
  $$ \left(p \geq 0; \,\, |z|<1,\,\, \Re(c)>\Re(b)>0  \right)$$
and
\begin{equation}\label{CHF-K}
 \Phi_{p,\nu}\,(b;\,c;\, z):=\sum_{n=0}^\infty\, \frac{B_{\nu}(b+n,\,c-b\,;\,p)}{B(b,\,c-b)}\,\frac{z^n}{n!}
\end{equation}
  $$ \left(p \geq 0; \,\,\Re(c)>\Re(b)>0  \right).$$

\vspace{3mm}

\noindent
{\bf Remark 5.}  The special case $\nu=0$ in (\ref{GHF-K}) and (\ref{CHF-K}) leads to the corresponding extensions
given in (\ref{GHF-p}) and (\ref{CHF-p}) by Chaudhry {\it et al.} \cite{Ch-Qa-Sr-Pa}.
\vspace{0.5cm}

\noindent 6.1.\ \ {\it Integral representations of the extended hypergeometric functions}
\vspace{0.3cm}

\noindent
In this section, we obtain integral representations of the extended Gauss hypergeometric and confluent hypergeometric functions as follows:
\vskip 3mm

\noindent
{\bf Theorem 7.} \emph{The following integral representations for the extended hypergeometric functions in}  (\ref{GHF-K}) and (\ref{CHF-K}) \emph{hold true}:
\[
F_{p,\nu}\,(a,\, b;\,c;\, z)=\sqrt{\frac{2p}{\pi}}\,\frac{1}{B(b,\,c-b)}\int_0^1\, t^{b-\frac{3}{2}}\,(1-t)^{c-b-\frac{3}{2}}\]
\begin{equation}\label{GHF-Integral-K}
   \times (1-zt)^{-a}\, K_{\nu+\frac{1}{2}}\left(\frac{p}{t(1-t)}\right)dt
\end{equation}
$$ \left(|\arg (1-z)|<\pi;\ p>0;\,\nu=0\,\,\mbox{and}\,\,p=0\,,\, \Re(c)>\Re(b)>0 \right)$$
\emph{and}
\[\Phi_{p,\nu}\,(b;\,c;\, z)= \sqrt{\frac{2p}{\pi}}\frac{1}{ B(b,\,c-b)}\int_0^1\, t^{b-\frac{3}{2}}\,(1-t)^{c-b-\frac{3}{2}}\,\exp \left(zt\right)\]
\begin{equation}\label{CHF-Integral-K}
   \times   K_{\nu+\frac{1}{2}}\left(\frac{p}{t(1-t)}\right)dt
\end{equation}
$$ \left(p>0;\,\nu=0\,\,\mbox{and}\,\,p=0\,,\Re(c)>\Re(b)>0 \right).$$

\noindent{\it Proof.}\ \ Substituting the definition of $B_{\nu}(x,\,y\,;\,p)$ in (\ref{Beta-function-K})
into (\ref{GHF-K}), we have
\[F_{p,\nu}\,(a,\, b;\,c;\, z)=\sqrt{\frac{2p}{\pi}}\,\frac{1}{B(b,\,c-b)}\int_{0}^{1}\, t^{b-\frac{3}{2}}\,(1-t)^{c-b-\frac{3}{2}}\hspace{4cm}\]
\begin{equation}\label{GHF-Integral-K-1}
\hspace{3cm}\times K_{\nu+\frac{1}{2}}\left(\frac{p}{t(1-t)}\right)\, \sum_{n=0}^\infty\, (a)_n\, \frac{(zt)^n}{n!}\,dt.
\end{equation}
Employing the binomial expansion
$$ (1-zt)^{-a} = \sum_{n=0}^\infty\, (a)_n\, \frac{(zt)^n}{n!}\qquad (|zt|<1)$$
in (\ref{GHF-Integral-K-1}), we obtain the integral in (\ref{GHF-Integral-K}).

A similar argument can be used to establish the integral representation of the extended confluent hypergeometric function in (\ref{CHF-Integral-K}) \hfill $\Box$.

\vspace{3mm}

\noindent
{\bf Remark 6.}  The special case $\nu=0$ of the integrals in (\ref{GHF-Integral-K}) and (\ref{CHF-Integral-K}) leads to the corresponding integral representations
given in Chaudhry {\it et al.} \cite{Ch-Qa-Sr-Pa}.

\vskip 3mm

\noindent
\textbf{Theorem 8.} \emph{The following representations of the extended hypergeometric functions for integer values} $\nu=n $ in (\ref{GHF-Integral-K}) and (\ref{CHF-Integral-K})
\emph{hold true}:
\begin{equation}\label{GHF-Integral-values}
F_{p,n}\,(a,\, b;\,c;\, z)= \sum_{m=0}^{n} \frac{(2p)^{-m}}{m!} \;\frac{\Gamma (n+m+1) }{\Gamma (n-m+1)}\frac{B(b+m,\,c-b+m)}{B(b,\,c-b)}F_{p}\,(a,\, b+m;\,c+2m;\, z)
\end{equation}
\emph{and}
\begin{equation}\label{CHF-Integral-values}
\Phi_{p,n}\,(b;\,c;\, z)= \sum_{m=0}^{n} \frac{(2p)^{-m}}{m!} \;\frac{\Gamma (n+m+1) }{\Gamma (n-m+1)}\frac{B(b+m,\,c-b+m)}{B(b,\,c-b)}\Phi_{p}\,(b+m;\,c+2m;\, z),
\end{equation}
\emph{where} $F_p$ \emph{and} $\Phi_p$ \emph{denote the extended functions defined in} (\ref{GHF-p}) \emph{and} (\ref{CHF-p}).
\vspace{0.2cm}

\noindent
{\it Proof.}\ \  Using the fact (\ref{K-formula}) in the integral representation (\ref{GHF-Integral-K}) and applying the definition (\ref{GHF-Integral-p}), we obtain the desired representation.

A similar argument can be used to establish the representation of the extended confluent hypergeometric function in (\ref{CHF-Integral-values}). \hfill $\Box$

\vspace{0.5cm}

\noindent 6.2.\ \ {\it Differentiation formulas}
\vspace{0.3cm}

\noindent
Differentiation formulas for the extended Gauss hypergeometric and confluent hypergeometric functions can be found by differentiating (\ref{GHF-K}) and (\ref{CHF-K}) with respect to $z$ as follows:

\vspace{3mm}

\noindent
{\bf Theorem 7.} \emph{The following differentiation formulas for the extended hypergeometric functions in} (\ref{GHF-K}) \emph{and}  (\ref{CHF-K})
\emph{hold true for} $p>0$  \emph{and non-negative integer} $n$:

\begin{equation}\label{GHF-Differentiation}
 \frac{d^n}{dz^n}\left\{F_{p,\nu}\,(a,\, b;\,c;\, z)\right\} =\frac{(a)_n\,(b)_n}{(c)_n}\,F_{p,\nu}\,(a+n,\, b+n;\,c+n;\, z)
  \end{equation}
\emph{and}
\begin{equation}\label{CHF-Differentiation}
 \frac{d^n}{dz^n}\left\{\Phi_{p,\nu}\,(b;\,c;\, z)\right\} =\frac{(b)_n}{(c)_n}\,\Phi_{p,\nu}\,(b+n;\,c+n;\, z)
\end{equation}

\noindent
{\it Proof.}\ \ Differentiating (\ref{GHF-K}) with respect to $z$, we have
\[
\frac{d}{dz}\left\{F_{p,\nu}\,(a,\, b;\,c;\, z)\right\} = \sum_{n=1}^\infty\,(a)_n\, \frac{B_{\nu}(b+n,\,c-b\,;\,p)}{B(b,\,c-b)}
  \, \frac{z^{n-1}}{(n-1)!}
\]
which, upon replacing $n$ by $n+1$ and using the facts  that
$$ B(b,\,c-b) = \frac{c}{b}\,B(b+1,\,c-b) \quad \mbox{and} \quad (\lambda)_{n+1}=\lambda\,(\lambda+1)_n,$$
yields
\[
 \frac{d}{dz}\left\{F_{p,\nu}\,(a,\, b;\,c;\, z)\right\} =\frac{a\,b}{c}\,F_{p,\nu}\,(a+1,\, b+1;\,c+1;\, z).
\]

The restriction $\Re(c)>\Re(b)>0$ on the above result may be removed by analytic continuation.
Repeated application of this process gives the
general form (\ref{GHF-Differentiation}).
A similar argument can be employed to establish (\ref{CHF-Differentiation}) for the extended confluent hypergeometric function. \hfill $\Box$

\vspace{3mm}

\noindent
{\bf Remark 8.}  The special case $\nu=0$ of (\ref{GHF-Differentiation}) and (\ref{CHF-Differentiation}) leads to the corresponding results given in Chaudhry {\it et al.} \cite{Ch-Qa-Sr-Pa}.

\vspace{0.5cm}

\noindent 6.3.\ \ {\it Transformation and summation formulas}
\vspace{0.3cm}

\noindent
In this section, we obtain  transformation and summation formulas for the extended Gauss hypergeometric and confluent hypergeometric functions as follows:

\vspace{3mm}

\noindent
{\bf Theorem 10.} \emph{The following transformation formulas for the extended hypergeometric functions hold true
when} $p>0$:

\begin{equation}\label{GHF-Transformation}
\hspace{1.5cm}F_{p,\nu}\,(a,\, b;\,c;\, z) =(1-z)^{-a}\,F_{p,\nu}\,\left(a,\, c-b;\,c;\, -\frac{z}{1-z}\right)\quad (|\arg (1-z)|<\pi)
\end{equation}
\emph{and}
\begin{equation}\label{CHF-Transformation}
\Phi_{p,\nu}\,(b;\,c;\, z)   = e^z\,\Phi_{p,\nu}\,(c-b;\,c;\,-z).
\end{equation}

\noindent
{\it Proof.}\ \ By writing
$$ \left[1-z(1-t) \right]^{-a} = (1-z)^{-a}\, \left(1 + \frac{z}{1-z}\,t  \right)^{-a}$$
and replacing $t$ by $1-t$ in (\ref{GHF-Integral-K}), we obtain
\[F_{p,\nu}\,(a,\, b;\,c;\, z) =  \frac{(1-z)^{-a}}{B(b,\,c-b)}\, \int_0^1\, t^{c-b-\frac{3}{2}}\,(1-t)^{b-\frac{3}{2}} \hspace{4cm}\]
\[\hspace{3cm}\times\left(1 + \frac{z}{1-z}\,t  \right)^{-a}\,K_{\nu+\frac{1}{2}}\left(\frac{p}{t(1-t)}\right)dt,\]
which establishes (\ref{GHF-Transformation}).

Again, replacing $t$ by $1-t$ in (\ref{CHF-Integral-K}), we find

\begin{equation}\label{CHF-Integral-K-1}
\Phi_{p,\nu}\,(b;\,c;\, z)   = \sqrt{\frac{2p}{\pi}}\frac{e^z}{ B(b,\,c-b)}\int_0^1\, t^{c-b-\frac{3}{2}}\,(1-t)^{b-\frac{3}{2}}\, \exp \left(-zt\right) K_{\nu+\frac{1}{2}}\left(\frac{p}{t(1-t)}\right)dt,
\end{equation}
which establishes (\ref{CHF-Transformation}). \hfill $\Box$

\vspace{3mm}

\noindent
{\bf Remark 9.} The special case $\nu=0$ in (\ref{GHF-Transformation}) and (\ref{CHF-Transformation}) leads to the corresponding extensions of Gauss' summation formula and Kummer's first transformation formula, respectively given in Chaudhry {\it et al.} \cite{Ch-Qa-Sr-Pa}.
\vspace{0.3cm}

\noindent
{\bf Theorem 11.} \emph{The following summation formula for the extended Gauss hypergeometric function holds true}:
\begin{equation}\label{GHF-summation}
F_{p,\nu}\,(a,\, b;\,c;\, 1)= \frac{B_{\nu}(b,\,c-a-b;\,p)}{B(b,\,c-b)}
\end{equation}
$$ \left(p>0;\,\nu=0\,\,\mbox{and}\,\,p=0,\,\,\Re(c-a-b)>0   \right).$$

\noindent
{\it Proof.}\ \ Setting $z=1$ in (\ref{GHF-Integral-K}) and using the definition (\ref{Beta-function-K}), we readily obtain the summation formula (\ref{GHF-summation}) \hfill $\Box$.
\vspace{3mm}

\noindent
{\bf Remark 10.} The formula (\ref{GHF-summation}) in the case $\nu=0$ and $p=0$ reduces to the well-known Gauss summation formula 
\[F(a,b;c;1)=\frac{\Gamma(c) \Gamma(c-a-b)}{\Gamma(c-a) \Gamma(c-b)}=\frac{B(b,c-a-b)}{B(b,c-b)} \qquad(\Re (c-a-b)>0).\]
\vspace{0.5cm}

\noindent 6.4.\ \ {\it A generating function for $F_{p,\nu}(a,b;c;z)$}
\vspace{0.3cm}

\noindent
We have the following theorem:
\vspace{3mm}

\noindent
{\bf Theorem 12.} \emph{The following generating function for $F_{p,\nu}(a,b;c;z)$ holds true}
\begin{equation}\label{e1}
\sum_{n=0}^\infty (a)_n\,F_{p,\nu}(a+n,b;c;z)\,\frac{t^n}{n!}=(1-t)^{-a} F_{p,\nu}\left(a,b;c;\frac{z}{1-t}\right).
\end{equation}
$$(p\geq 0,\,\,|t|<1).$$

\noindent
{\it Proof.}\ \ Let the left-hand side of (\ref{e1}) be denoted by $S$.Then, from the definition (\ref{GHF-K}) we have
\begin{eqnarray*}
S&=&\sum_{n=0}^\infty (a)_n
\left\{\sum_{k=0}^\infty (a+n)_k\,\frac{B_\nu(b+k,c-b;p)}{B(b,c-b)}\,\frac{z^k}{k!}\right\} \frac{t^n}{n!}\\
&=&\sum_{k=0}^\infty (a)_k\,\frac{B_\nu(b+k,c-b;p)}{B(b,c-b)} \left\{\sum_{n=0}^\infty (a+k)_n\,\frac{t^n}{n!}\right\}\frac{z^k}{k!}
\end{eqnarray*}
upon reversal of the order of summation and use of the identity $(a)_n(a+n)_k=(a)_k(a+k)_n$.

Since by the binomial theorem
\[(1-t)^{-a-k}=\sum_{n=0}^\infty (a+k)_n\,\frac{t^n}{n!}\qquad (|t|<1),\]
identification of the series over $k$ from (\ref{GHF-K}) as $F_{p,\nu}(a,b;c;z/(1-t))$ then leads to the assertion 
in (\ref{e1}). \hfill $\Box$

\vspace{0.5cm}

\noindent 6.5.\ \ {\it Asymptotic behaviour for large $z$}
\vspace{0.3cm}

\noindent
We present the asymptotic behaviour of $\Phi_{p,\nu}(b;c;z)$ for $|z|\rightarrow\infty$ in $\Re(z)>0$ when it is supposed that $p>0$. The major contribution to the integrand in (\ref{CHF-Integral-K-1}) for large $|z|$ arises from the neighbourhood of $t=0$, where the Bessel function may be approximated by its leading asymptotic behaviour
\[K_{\nu+\frac{1}{2}}\left(\frac{p}{t(1-t)}\right)\sim\sqrt{\frac{\pi}{2p}}\, t^\frac{1}{2}(1-t)^\frac{1}{2} \exp\left(-\frac{p}{t}-\frac{p}{1-t}\right) \qquad (t\rightarrow 0).\]
Then, from (\ref{CHF-Integral-K-1}), we obtain with $t=z\tau$
\begin{eqnarray}\Phi_{p,\nu}(b;c;z)&\simeq& \frac{z^{b-c} e^{z-p}}{B(b,c-b)} \int_0^z \tau^{c-b-1}\exp(-\tau-\frac{p}{\tau})\,d\tau\nonumber\\
&=&\frac{z^{b-c} e^{z-p}}{B(c,c-b)} \,\gamma(c-b,z;pz)  \qquad (|z|\rightarrow\infty,\ \Re(z)>0,\ p>0),
\end{eqnarray}
where $\gamma(\alpha,x;p)$ is the extended incomplete gamma function defined in (\ref{Incomplete-gamma-p}). We note that the leading term is independent of the parameter $\nu$ and that the factor $e^{-p}$ was omitted in the case $\nu=0$ in \cite[Eq.~(8.7)]{Ch-Qa-Sr-Pa}.
In the limit $p\rightarrow 0$, we recover the well-known asymptotic behaviour of the ordinary confluent hypergeometric function
\[\Phi(b;c;z) \sim z^{b-c} e^z\,\frac{\gamma(c-b,z)}{B(b,c-b)}\sim\frac{\Gamma(c)}{\Gamma(b)}\,z^{b-c} e^z
\qquad (|z|\rightarrow\infty,\ \Re(z)>0).\]

For the function $F_{p,\nu}(a,b;c;z)$, we apply a slight modification of the analysis described in \cite[\S 16.3]{WW} to the integral (\ref{GHF-Integral-K}), where for simplicity in presentation we shall suppose the parameters $a$, $b$, $c$ and $\nu$  to be real. With $\xi=-z$, we write 
\[F_{p,\nu}(a,b;c;z)=\sqrt{\frac{2p}{\pi}}\,\frac{(-z)^{-a}}{B(b,c-b)} \int_0^1t^{b-a-\frac{3}{2}}(1-t)^{c-b-\frac{3}{2}}\left(1+\frac{1}{\xi t}\right)^{\!\!-a}K_{\nu+\frac{1}{2}}\left(\frac{p}{t(1-t)}\right)dt,\] 
where we employ the series expansion for positive integer $n$
\[\left(1+\frac{1}{\xi t}\right)^{\!\!-a}=\sum_{k=0}^{n-1}\frac{(-1)^k(a)_k}{(\xi t)^{k}k!}+R_n(t,\xi)\]
with
\[R_n(t,\xi)=\frac{(-1)^n(a)_n}{(n-1)!} \left(1+\frac{1}{\xi t}\right)^{\!\!-a} \int_0^{1/(\xi t)} u^{n-1}(1+u)^{a-1}du.\]
The finite series yields the contribution to $F_{p,\nu}(a,b;c;z)$ given by
\[\sqrt{\frac{2p}{\pi}}\,\frac{(-z)^{-a}}{B(b,c-b)}\sum_{k=0}^{n-1} \frac{(a)_k}{k!} z^{-k} \int_0^1 t^{b-a-k-\frac{3}{2}}(1-t)^{c-b-\frac{3}{2}} K_{\nu+\frac{1}{2}}\left(\frac{p}{t(1-t)}\right)dt\]
\[=\frac{(-z)^{-a}}{B(b,c-b)} \sum_{k=0}^{n-1} \frac{(a)_k}{k!} z^{-k} B_\nu(b-a-k,c-b;p).\]

Now
\begin{eqnarray*}
|R_n(t,\xi)|&\leq&\frac{|(a)_n|}{(n-1)!} \left|\left(1+\frac{1}{\xi t}\right)^{\!\!-a}\right|\int_0^{1/(t|z|)} u^{n-1}(1+u)^{|a|}du\\
&\leq&\frac{|(a)_n|}{(n-1)!} \frac{\left(1+t^{-1}\right)^{|a|}}{(|z|t)^n(\sin \epsilon)^{|a|}} \int_0^1 w^{n-1}\left(1+\frac{w}{\xi t}\right)^{\!\!|a|} dw \\
&\leq& \frac{|(a)_n|}{n!} \frac{(1+t^{-1})^{2|a|}}{(|z|t)^n (\sin \epsilon)^{|a|}},
\end{eqnarray*}
where, with $|z|>1$, $|\arg\,\xi|\leq\pi-\epsilon$, we have used the bounds
\[1\leq |1+(\xi t)^{-1}|\leq 1+t^{-1}\ \ (\Re(\xi)\geq 0),\qquad |1+(\xi t)^{-1}|\geq \sin \epsilon\ \ (\Re(\xi)\leq 0).\]
It then follows that, when $|z|>1$,
\[\left|\int_0^1 t^{b-a-\frac{3}{2}}(1-t)^{c-b-\frac{3}{2}} R_n(t,\xi)K_{\nu+\frac{1}{2}}\left(\frac{p}{t(1-t)}\right)dt\right|\hspace{5cm}\]
\[\leq \frac{|(a)_n|}{n! |z|^n (\sin \epsilon)^{|a|}}\int_0^1 t^{b-a-n-\frac{3}{2}}(1-t)^{c-b-\frac{3}{2}}(1+t^{-1})^{2|a|}K_{\nu+\frac{1}{2}}\left(\frac{p}{t(1-t)}\right)dt
=O(z^{-n}),\]
since the integral is convergent when $p>0$ and is independent of $z$. The constant implied in the $O$-symbol becomes infinite as $\epsilon\rightarrow 0$.

Hence, we obtain the following theorem:
\vspace{3mm}

\noindent
{\bf Theorem 13.} \emph{For positive integer $n$ and $p>0$, we have the asymptotic expansion} 
\begin{equation}
F_{p,\nu}(a,b;c;z)=\frac{(-z)^{-a}}{B(b,c-b)}\sum_{k=0}^{n-1}\frac{(a)_k}{k!} z^{-k} B_\nu(b-a-k,c-b;p)+O(z^{-n})
\end{equation}
\emph{for} $|z|\rightarrow\infty$ \emph{in} $|\arg (-z)|<\pi$.

\vspace{0.6cm}

\begin{center}
{\bf 7. \ Mellin Transforms Representations} 
\end{center}
\setcounter{section}{7}
\setcounter{equation}{0}
\renewcommand{\theequation}{\arabic{section}.\arabic{equation}}
\noindent
In this section, we obtain Mellin transform representations for the extended Gauss hypergeometric and confluent hypergeometric functions as follows:
\vspace{3mm}

\noindent
{\bf Theorem 14.} \emph{The following Mellin transformation formulas for the extended hypergeometric functions in} (\ref{GHF-K}) and (\ref{CHF-K})
\emph{hold true}:
\[{\mathcal M}\left\{F_{p,\nu}\,(a,\, b;\,c;\, z):p \to s\right\}\hspace{9cm}\]
\begin{equation}\label{Mellin-transform-GHF} =\frac{2^{s-1}}{\surd\pi}\,\Gamma(\frac{s-\nu}{2})\Gamma(\frac{s+\nu+1}{2})\;\frac{B(b+s,\,c+s-b)}{B(b,\,c-b)}\,F(a,\, b+s;\,c+2s;\, z)
\end{equation}
\emph{and}
\[{\mathcal M}\left\{\Phi_{p,\nu}\,(b;\,c;\, z):p \to s\right\}\hspace{9cm}\]
\begin{equation}\label{Mellin-transform-CHF} =\frac{2^{s-1}}{\surd\pi}\,\Gamma(\frac{s-\nu}{2})\Gamma(\frac{s+\nu+1}{2})\;\frac{B(b+s,\,c+s-b)}{B(b,\,c-b)}\,\Phi(b+s;\,c+2s;\, z)
\end{equation}
$$(\Re(s-\nu)>0,\ \Re(s+\nu)>-1),$$
provided that each member of (\ref{Mellin-transform-GHF}) and (\ref{Mellin-transform-CHF}) exists.

\noindent
{\it Proof.} \ Using the definition (\ref{Mellin-transform}) of the
Mellin transform, we find from (\ref{GHF-Integral-K}) that
\[{\mathcal M}\left\{F_{p,\nu}\,(a,\, b;\,c;\, z):p \to s\right\}\hspace{8cm}\]
\[=\int_0^{\infty}p^{s-1}
\left\{\sqrt{\frac{2p}{\pi}}\,\frac{1}{B(b,\,c-b)}\int_0^1\, t^{b-\frac{3}{2}}\,(1-t)^{c-b-\frac{3}{2}}(1-zt)^{-a}\,   K_{\nu+\frac{1}{2}}\left(\frac{p}{t(1-t)}\right)\,dt\right\}dp\]
\[=\sqrt{\frac{2}{\pi}}\,\frac{1}{B(b,\,c-b)}\int_0^{\infty}p^{s-\frac{1}{2}}\left\{
\int_0^1  t^{b-\frac{3}{2}}\,(1-t)^{c-b-\frac{3}{2}}(1-zt)^{-a}\,K_{\nu+\frac{1}{2}}\left(\frac{p}{t(1-t)}\right)  dt\right\}dp.\]
Then, upon interchanging the order of integration by absolute convergence of the integrals
involved, we obtain
\vspace{2mm}
\[{\mathcal M}\left\{F_{p,\nu}\,(a,\, b;\,c;\, z):p \to s\right\}\hspace{9cm}\]
\[ =\sqrt{\frac{2}{\pi}}\,\frac{1}{B(b,\,c-b)}\int_0^1t^{b-\frac{3}{2}}\,(1-t)^{c-b-\frac{3}{2}}(1-zt)^{-a}
\left\{\int_0^{\infty}p^{s-\frac{1}{2}}K_{\nu+\frac{1}{2}}\left(\frac{p}{t(1-t)}\right)dp\right\} dt\]
\[=\int_0^1t^{b+s-1}\,(1-t)^{c-b+s-1}(1-zt)^{-a}
\left\{\sqrt{\frac{2}{\pi}}\,\int_0^{\infty}w^{s-\frac{1}{2}}K_{\nu+\frac{1}{2}}\left(w\right) dw\right\} dt,\]
where we have set $w=p t/(1-t)$ in the inner $p$-integral.  Then, upon use of (\ref{K-integral}) to evaluate the inner integral followed by (\ref{GHF-Integral-p}) (with $p=0$),
we finally obtain (\ref{Mellin-transform-GHF}) subject to the conditions stated. This completes
our derivation of the Mellin transform Formula asserted by Theorem 12.

A similar argument can be employed to establish (\ref{Mellin-transform-CHF}) for the extended confluent hypergeometric function.
\hfill $\Box$

\end{document}